\documentclass[journal,twoside,letterpaper]{ieeeconf}

\usepackage{amssymb}
\usepackage{bm}
\usepackage{amsmath}
\usepackage{graphicx}
\usepackage{xcolor}
\usepackage{hyperref}
 \hypersetup{colorlinks=true, linkcolor=blue, breaklinks=true, urlcolor=blue, citecolor=blue}
\usepackage{mathtools}

\newcommand{\Id}{{\rm Id}}
\newcommand{\JFF}{{\bm J}_{\cal FF}}
\newcommand{\JFS}{{\bm J}_{\cal FS}}
\newcommand{\Jred}{{\bm J}_{\rm red}}
\newcommand{\JSF}{{\bm J}_{\cal SF}}
\newcommand{\JSS}{{\bm J}_{\cal SS}}

\begin{document}

\title{Nontrivial Kron Reduction for Power Grid Dynamics Modeling}
\author{Laurent Pagnier, \emph{IEEE Member}, Robin Delabays, and Melvyn Tyloo
\thanks{L. Pagnier is with the Department of Mathematics, University of Arizona, Tucson, AZ 85721, United States of America.}
\thanks{R. Delabays is with the Institute of Sustainable Energy, School of Engineering, University of Applied Sciences and Arts of Western Switzerland (e-mail: robin.delabays@hevs.ch).}
\thanks{M. Tyloo is with the Theoretical Division and Center for Nonlinear Studies (CNLS), Los Alamos National Laboratory, Los Alamos, NM 87545, USA (e-mail: melvyn.tyloo@gmail.com).}
\thanks{RD was supported by the Swiss National Science Foundation under grant no.  	200021\_215336.}
\thanks{MT was supported by the Laboratory Directed Research and Development program of Los Alamos National Laboratory under project numbers 20220797PRD2 and 20220774ER and by U.S. DOE/OE as part of the DOE Advanced Sensor and Data Analytics Program.}
\thanks{The authors thank Andrey Lokhov for useful discussions.}
}

\maketitle
\thispagestyle{empty}

\begin{abstract}
 The Kron reduction is used in power grid modeling when the analysis can -- supposedly -- be restricted to a subset of nodes.  Typically, when one is interested in the phases' dynamics, it is common to reduce the load buses and focus on the generators' behavior. The rationale behind this reduction is that voltage phases at  load buses adapt quickly to their neighbors' phases and, at the timescale of generators, they have virtually no dynamics. We show that the dynamics of the Kron-reduced part of a network can have a significant impact on the dynamics of the non-reduced buses. Therefore, Kron reduction should be used with care and, depending on the context, reduced nodes cannot be simply ignored.  We demonstrate that the noise in the reduced part can unexpectedly affect the non-reduced part, even under the assumption that nodal disturbances are independent. Therefore, the common assumption that the noise in the non-reduced part is uncorrelated may lead to inaccurate assessments of the grid's behavior.  To cope with such shortcomings of the Kron reduction, we show how to properly incorporate the contribution of the reduced buses into the reduced model using the Mori-Zwanzig formalism.
\end{abstract}


\maketitle

\section{Introduction}
The resilience of networked systems is of primal importance to ensure their reliability and operations.
It is especially the case for electric power grids, which have attracted tremendous attention over the last decade due to the ongoing energy transition~\cite{brummitt2013transdisciplinary,smith2022effect}. 
One way to evaluate the resilience of a system is by quantifying its response to an external input~\cite{bamieh2012coherence,siami2014systemic,tegling2015price,grunberg2018performance}. 
In a high-voltage transmission grid, for instance, one can wonder about the voltage phase response to a noisy injection of power coming from a renewable energy source. 

Such transmission grids are usually large-scale systems.
As a convenient way to render their analysis more tractable, it has become standard to first apply a \textit{Kron reduction}~\cite{kron1939tensor, dorfler2013kron}. 
By taking the Schur complement~\cite[Sec.~0.8.5]{Hor94} of the dynamics matrix of a system, the latter transformation only retains the synchronous generators while reducing the loads into effective line susceptances. 
Then, based on this initial reduction, previous studies focused on the resilience to disturbances occurring at the non-reduced buses. 

The rationale underlying the validity of the Kron reduction boils down to that of timescale separation. 
Indeed, loads typically respond much faster than conventional plants to frequency disturbances. 
Loads then appear as {\it passive nodes} in the timescale of the plants. 

While the analysis of the Kron-reduced system allows to approximate the vulnerability of synchronous generators, we claim that such an approach is insufficient to accurately and fully evaluate the resilience to external inputs. 
Obviously, the reduced buses, which are typically loads, are subject to disturbances too. 
In fact, their dynamics is mostly fluctuating in time, as they are dictated by the consumption of power. 
We show in this manuscript that, although the noise at the reduced nodes is spatially uncorrelated, (i) its impact on non-reduced nodes may be correlated and (ii) the magnitude of said impact may be hard to predict. 
Importantly, we show that summarizing the impact of the reduced nodes' noise on the non-reduced nodes to a white noise can be detrimental to an accurate resilience estimation.

We propose a unified framework to assess the impact of disturbances both at the reduced and non-reduced buses. 
First, using Mori-Zwanzig formalism~\cite{mori1965transport,zwanzig1973nonlinear,Tyl24}\,, we take an unconventional approach to account for the timescale separation assumed by the Kron reduction. We derive a general expression for the time-evolution of the slow buses as a function of the fast ones. 
By expanding over the timescale separation parameters ($\epsilon$ below), this naturally leads to the reduction of the grid dynamics to degrees of freedom corresponding to the synchronous generators. 
Applying this formalism, we also account for the reduction of the disturbance, and explicitly calculate their effect on the non-reduced buses. 

Second, we assess the resilience of the grid by evaluating the variance of the frequency deviations induced by time-correlated noisy inputs. 
In particular, we can compare the contribution from reduced and non-reduced buses to the variance of the frequency deviation.
We show, through analytical and numerical examples how misleading simulation of Kron-reduced systems can be, if disturbances in the reduced part are neglected. 

\section{Preliminaries}\label{sec2}
Throughout the manuscript, we use $\bm{1}_N$ to denote the $N$-dimensional vector of ones and $\Id_N$ for the identity matrix of size $N$. In the following, we sometimes refer to buses as nodes.

We model the voltage phase dynamics with the structure preserving model~\cite{Ber81} of the swing equations~\cite[Sec.~3.9.2]{Kun22}. 
Neglecting line conductances and voltage variations as a first approximation, the dynamics of voltage phases resort to that of a network of $N$ nonlinearly coupled phases $\theta_i\in(-\pi,\pi]$, 
\begin{align}\label{eq1}
\begin{split}
    m_i\,\ddot{\theta}_i+d_i\,\dot{\theta}_i = p_i - \sum_{j=1}^N b_{ij}\,\sin(\theta_i-\theta_j) + \eta_i\,,
\end{split}
\end{align}
for $i=1,...\,N$\,. 
The effective inertia and damping at each bus $i$ are summarized in the constants $m_i$ and $d_i$ respectively, and $p_i$ denotes the power injected or withdrawn at that bus. The coupling between the buses is given by the elements of the adjacency matrix $b_{ij}=B_{ij}|V_i||V_j|$ where $B_{ij}$ is the susceptance of the transmission line connecting buses $i$ and $j$ and $|V_i|$ is the amplitude (assumed to be constant) of the voltage at bus $i$\,. 
The noise induced by local power variations $\eta_i$ is modelled by a time-correlated space-uncorrelated noise, 
\begin{align}
 \langle \eta_i(t)\eta_j(t') \rangle &= \sigma_{i}^2\delta_{ij}\exp[-|t-t'|/\tau_i]\, ,
\end{align}
where $\tau_i$ is the typical correlation time of the noise at the $i$th node and $\delta_{ij}$ is the Kronecker symbol. 
Note that one could also consider space-correlated noise~\cite{ronellenfitsch2018optimal,hindes2023stability}. 
Here, for the sake of clarity, we stick to the common assumption of space-uncorrelated noise and demonstrate how it is affected by the reduction.

\section{Timescale Separation in the Swing Equations}\label{sec3}
Timescale separation is the rationale underlying the use of the Kron reduction. 
We formalize here what is meant by timescale separation. 
Then, we present the Mori-Zwanzig formalism that allows us to deal with the nontrivial impact of reduced nodes. 

\subsection{Timescale Separation and Reduced Dynamics}
In view of introducing the Kron reduction of the system, we assume that we have two sets of buses that we denote $\mathcal{F}$ and $\mathcal{S}$\,, respectively with $N_{\mathcal{F}}$ and $N_{\mathcal{S}}$ nodes.
Time scale separation can be summarized in terms of inertia and damping properties, 
\begin{eqnarray}
d_i,m_i\propto \begin{cases}
\overline{d},\overline{m} & i\in\mathcal{F}\, ,\\
\underline{d},\underline{m} & i\in\mathcal{S}\, ,
\end{cases}
\end{eqnarray}
with $\overline{d}\ll \underline{d}$\,, $\overline{m}\ll \underline{m}$\,. 
The latter means that buses belonging to $\mathcal{S}$ have a much slower intrinsic timescale than those belonging to $\mathcal{F}$\,. 
The buses in ${\cal F}$ are the ones that will be reduced by Kron reduction.

In the following, we focus on the dynamics of the buses in the slow component. 
Within the assumption of timescale separation, adapting parameter definitions accordingly, one can rewrite~(\ref{eq1}) as
\begin{align}\label{eq2}
\!\!\begin{split}
    m_i\,\ddot{\theta}_i+d_i\,\dot{\theta}_i &= p_i - \sum_{j=1}^N b_{ij}\,\sin(\theta_i-\theta_j) + \eta_i\, , \,i\in\mathcal{S},\\
   \epsilon\,(m_i\,\ddot{\theta}_i+d_i\,\dot{\theta}_i) &= p_i - \!\!\sum_{j=1}^N b_{ij}\,\sin(\theta_i-\theta_j) + \eta_i\, , \,i\in\mathcal{F},
\end{split}
\end{align}
where we defined $\underline{m}/\overline{m} = \underline{d}/\overline{d}=\epsilon^{-1}$\,.
In the limit $\epsilon\rightarrow 0$\,, the buses within $\mathcal{F}$ instantaneously adapt their phases. In the following we show how to reduce these buses. 

\subsection{System Response and Mori-Zwanzig Analysis}
In order to assess the resilience of the system, we want to analyze the frequency deviation from nominal value when subject to noise. 
We therefore consider the system~\eqref{eq2} in the vicinity of a stable fixed point $\bm{\theta}^*$.
In particular, we are interested in the time-evolution of the phase deviations $x_i(t)=\theta_i(t)-\theta_i^{*}$ for $i\in\mathcal{S}$ and $y_i(t)=\theta_i(t)-\theta_i^{*}$ for $i\in\mathcal{F}$ whose dynamics at the first order read,
\begin{align}\label{eq4}
{\bm M}\begin{bmatrix}
\ddot{\bf x}\\
\ddot{\bf y}
\end{bmatrix}+{\bm D}\begin{bmatrix}
\dot{\bf x}\\
\dot{\bf y}
\end{bmatrix}=
\begin{bmatrix}
\JSS & \JSF \\
\JFS & \JFF
\end{bmatrix}
\begin{bmatrix}
{\bf x}\\
{\bf y}
\end{bmatrix}+
\begin{bmatrix}
{\bm \eta}_\mathcal{S}\\
{\bm \eta}_\mathcal{F}
\end{bmatrix}
\end{align}
where we defined the diagonal inertia and damping matrices,
\begin{align}\label{defMD}
{\bm M} &= \begin{bmatrix}
 {\bm M}_{\mathcal{S}} & {\bm 0} \\
 {\bm 0} & \epsilon\,{\bm M}_{\mathcal{F}}
\end{bmatrix}\, , & 
{\bm D} &= \begin{bmatrix}
 {\bm D}_{\mathcal{S}} & {\bm 0} \\
 {\bm 0} & \epsilon\,{\bm D}_{\mathcal{F}}
\end{bmatrix}\, ,
\end{align}
and the Jacobian matrix of the system
\begin{eqnarray}\label{eqjac}
J_{ij}=\begin{cases}
 b_{ij}\cos(\theta_i^*-\theta_j^*)& i\neq j\\
 -\sum_{k=1}^Nb_{ik}\cos(\theta_i^*-\theta_k^*)& i=j\,,
\end{cases}
\end{eqnarray}
which is a Laplacian matrix when phase differences are between $-\frac{\pi}{2}$ and $\frac{\pi}{2}$\,. 
Notice that for undirected coupling as we consider in the following, one has $\JFS = \JSF^\top$\,. 

Ultimately, one would like to obtain a closed-form expression for the frequency at the generator buses as these are the state variables relevant for the stability of the grid. To obtain such expression, one can use Mori-Zwanzig formalism~\cite{mori1965transport,zwanzig1973nonlinear,Tyl24} which was developed in statistical mechanics to derive the time-evolution of relevant variables while treating the rest of the variables as input noise. One can follow a similar approach with $\bf x$ and $\bf y$ being respectively the relevant and irrelevant variables, see \cite{Tyl24} for more details. 

For the sake of simplicity, we assume that the inertia and damping parameters are homogeneous, i.e., $m_i=m$\,, $d_i=d$ $\forall i$\,. 
Nevertheless, simulations shown in Sec.~\ref{sec:sim} demonstrate that similar results hold for heterogeneous parameters.
One can express the first row of~(\ref{eq4}) as,
\begin{align}\label{eqmz}
    &m\,\ddot{x}_i+d\,\dot{x}_i = \sum_{j=1}^{N_{\cal S}}(\JSS)_{ij} x_j + \left(\eta_{\cal S}\right)_i + \sum_{j=1}^{N_{\cal F}}\left(\JSF\right)_{ij}y_j\, ,
\end{align}
with
\begin{align}\label{eqy}
    {y_i}(t) = &m^{-1}\epsilon^{-1}\sum_{\alpha=1}^{N_{\mathcal{F}}}e^{-(\gamma+\Gamma_\alpha)t/2}\int_0^t e^{\Gamma_\alpha t_1}\\
    &\times\int_0^{t_1}[\JFS{\bf x}(t_2) + {\bm \eta}_{\mathcal{F}}]\cdot {\bf w}_\alpha e^{(\gamma-\Gamma_\alpha)t_2/2} {\rm d}t_2 {\rm d}t_1 {w}_{\alpha,i}\,,\nonumber
\end{align}
where $\Gamma_\alpha=\sqrt{\gamma^2-4\nu_\alpha/(m\epsilon)}$ and we denoted ${\bf w}_{\alpha}$ the eigenvectors of $\JFF$ with corresponding eigenvalues $\nu_\alpha<0$\,. 
Equation~(\ref{eqmz}) is the equation of motion of the relevant variables.
The expressions (\ref{eqmz}) and (\ref{eqy}) provide a general framework to analyse how the irrelevant variables impact the relevant ones. 
In the particular case where the two sets of variables evolve on different timescales, one can expand~(\ref{eqmz}) and~(\ref{eqy}) over the timescale parameter $\epsilon$ to perturbatively calculate the effect of $\bf y$ on $\bf x$\,. 
Here we are interested in the limit $\epsilon\rightarrow 0$ and therefore only consider the non-vanishing leading order in $\epsilon$\,, which is the zeroth order. Note however that one could calculate the higher order correction by a Taylor expansion in $\epsilon$\,. 
This limit induces Dirac-delta distributions between $t_1$ and $t$ and between $t_2$ and $t_1$ in~(\ref{eqy})~\cite{Tyl24}. 
Indeed, one has,
\begin{align}\label{eqylim}
    \lim_{\epsilon\rightarrow 0}m^{-1/2}\epsilon^{-1/2}e^{-(\gamma+\Gamma_\alpha)(t-t')/2} = {\nu_\alpha}^{-1/2}\delta(t-t')\,.
\end{align}
Therefore, noticing that $\JFF^{-1}=\sum_{\alpha=1}^{N_\mathcal{F}}{\nu_\alpha}^{-1}{\bf w}_{\alpha}{\bf w}_{\alpha}^\top$\,, one can rewrite~(\ref{eqmz}) in a matrix-vector form as,
\begin{align}\label{eq5}
    &{\bm M}\,\ddot{\bf x} + {\bm D}\,\dot{\bf x} \nonumber \\
    &=\JSS\,{\bf x} - \JSF\, \JFF^{-1}\,\JFS\,{\bf x} +  {\bm \eta}_{\mathcal{S}}- \JSF\, \JFF^{-1}\,{\bm \eta}_{\mathcal{F}} \\
    &=: \Jred\,{\bf x} + {\bm \xi}  \,,\nonumber
\end{align}
where in the last line we defined the reduced Jacobian $\Jred := \JSS - \JSF\, \JFF^{-1}\,\JFS$\,, and denoted the noise term as ${\bm \xi}:={\bm \eta}_{\mathcal{S}}- \JSF\, \JFF^{-1}\,{\bm \eta}_{\mathcal{F}}$\,. 
The reduced Jacobian is still a Laplacian matrix~\cite{dorfler2013kron}.
The dynamics of the non-reduced buses is then governed by~(\ref{eq5}) where the effective noise at the $i$th bus is a combination of the noise at the $i$th bus and of a superposition of the noise inputs at buses belonging to the fast component. 
Therefore, in general $\bm \xi$ is correlated over the different synchronous generators.

Now that we have performed the reduction using the Mori-Zwanzig formalism, let's compare \eqref{eq5} with the conventional Kron-reduced model. The reduced Jacobian, $\bm{J}_{\rm red}$, is given by the same expression in both cases. However, with the present formalism, the validity of this approximation is well-defined by the underlying assumptions. Since the Kron reduction is inherently a static analysis method, it does not provide a clear approach to properly define nodal noise signals. As a result, uncorrelated inputs, $\bm{\eta}_\mathcal{S}$, are often used. These differ from the noise signal $\bm{\xi}$ that we have derived in general.

The linear system~(\ref{eq5}) can be solved by expanding over the eigenmodes of $\Jred$ denoted ${\bf u}_{\alpha}$\,, with corresponding eigenvalues $\lambda_\alpha$\,, $\alpha=1,...N_{\mathcal{S}}=|\mathcal{S}|$\,. 
As $\Jred$ is also the negative of a Laplacian matrix, one has that $0=\lambda_1\ge ... \ge \lambda_{N_{\mathcal{S}}}$ with $u_{1,i}=1/\sqrt{N_{\mathcal{S}}}$\,. 
The general solution to~(\ref{eq5}) is given by,
\begin{align}\label{eqsol}
x_i(t) &= \sum_{\alpha=2}^{N_{\mathcal{S}}}(m\Gamma_\alpha)^{-1}e^{-(\gamma+\Gamma_\alpha)t/2}\\
&\times\int_{0}^{t}\left[e^{(\gamma+\Gamma_\alpha)t'/2} + e^{(\gamma-\Gamma_\alpha)t'/2}\right]{\bf u}_{\alpha}\cdot{\bm \xi}(t')\,{\rm d}t'u_{\alpha,i}\,,\nonumber
\end{align}
for $i=1,...,N_{\mathcal{S}}$ and where we omitted the first mode in the sum as any perturbation along it does not modify the system. In the following we only consider the dynamics orthogonal to the first mode. Eventually, we are interested in the variance of the frequency deviations. Using the expression in~(\ref{eqsol}) and taking its time-derivative, one can calculate the moments of the frequency deviations.

\section{Resilience of the grid}
We are now ready to derive the correction term to the non-reduced nodes' dynamics after Kron reduction. 
Furthermore, we provide a couple of idealized examples where the impact of the Kron reduction is rather intuitive. 

\subsection{Fluctuations From the Synchronized State}

Various characteristics of the response can be used to determine the resilience of the coupled nodes. When subject to stochastic inputs, a natural choice is to evaluate the magnitude of the deviations from the synchronized fixed point by calculating the variance of the frequency deviations. As discussed in Sec.~\ref{sec2}, we consider time-correlated noise inputs with different typical correlation times in each component, denoted as $\tau_{\mathcal{F}}$ and $\tau_{\mathcal{S}}$ respectively. 

The \textit{center of inertia} variance of the frequency deviations orthogonal to the zero mode ${\bf u}_1$ in the slow component is calculated from~(\ref{eqsol}) and reads, in the long-time limit, 
\begin{align}\label{eqfr}
     &\big\langle \dot x_i^2\big\rangle_{\rm COI} 
     \coloneqq \Big\langle \big(\dot x_i - N_\mathcal{S}^{-1}\sum_k^{N_\mathcal{S}} \dot x_k\big)^2\Big\rangle \\
     &= \sum_{\alpha,\beta=2}^{N_{\mathcal {S}}}\sum_{j=1}^{N_{\mathcal {S}}}\Big[U_{ij}^{\alpha\beta}\mathcal{H}(\lambda_\alpha, \lambda_\beta, \tau_\mathcal{S}, \gamma) 
     + \overline{\Gamma}^{\alpha\beta}_i \mathcal{H}(\lambda_\alpha, \lambda_\beta, \tau_\mathcal{F}, \gamma)\Big]\,,\nonumber  
\end{align}
where we defined the shorthand notations
\begin{align}
\begin{split}
    {U_{ij}^{\alpha\beta}} &= {\sigma_{\mathcal{S}}}_j^2u_{\alpha,i}u_{\alpha,j}u_{\beta,i}u_{\beta,j}\, , \\
    \overline{\Gamma}^{\alpha\beta}_i &= \Gamma_{\alpha\beta}\,u_{\alpha,i}u_{\beta,i}\, , \\
    \Gamma_{\alpha\beta} &= {\bf u}^\top_{\alpha} \JSF \JFF^{-1}\, {\rm diag}\big[\bm \sigma_{\mathcal{F}}^2\big]\, \JFF^{-1} \JFS {\bf u}_{\beta}\, ,
\end{split}
\end{align}
and 
\begin{align}
    \mathcal{H}(\lambda_\alpha, \lambda_\beta, \tau, \gamma) &= \nonumber \\
    &\hspace{-2.5cm}\frac{ \frac{1}{2}\left[2 \gamma m \overline{\Lambda}_{\alpha\beta}( \gamma  \tau +1)+4 \gamma \lambda_\alpha \lambda_\beta \tau^2-\tau \underline{\Lambda}_{\alpha\beta}^2\right]}{ \left(2 \gamma^2 m \overline{\Lambda}_{\alpha\beta}+\underline{\Lambda}_{\alpha\beta}^2\right) \left(\gamma m \tau+\lambda_\alpha \tau^2+m\right) \left(\gamma m \tau+\lambda_\beta \tau^2+m\right)}, \nonumber
\end{align}
with
\begin{align}
    \overline{\Lambda}_{\alpha\beta} &= (\lambda_\alpha+\lambda_\beta)\, , &
    \underline{\Lambda}_{\alpha\beta} &= (\lambda_\alpha-\lambda_\beta)\, .
\end{align}

\begin{figure}[h!]
    \centering
    \includegraphics[width=0.6\linewidth]{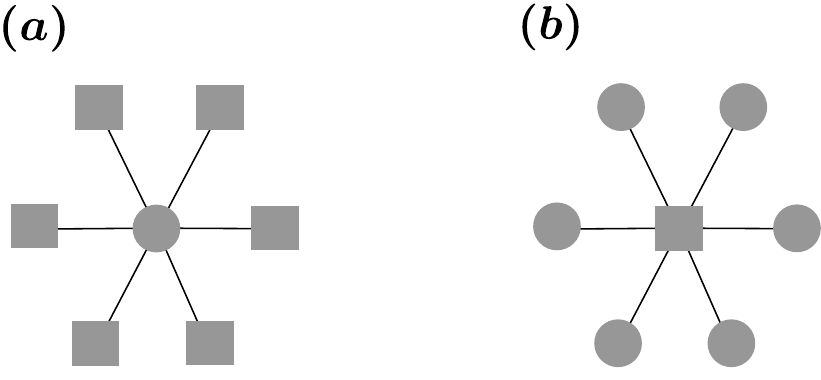}
    \caption{Star network configurations showing (a) a load surrounded by generators and (b) a generator surrounded by loads. Circles represent loads (fast buses), and squares represent generators (slow buses).}
    \label{fignet}
\end{figure}
\begin{figure*}[h!]
    \centering
    \includegraphics[width=0.90\linewidth]{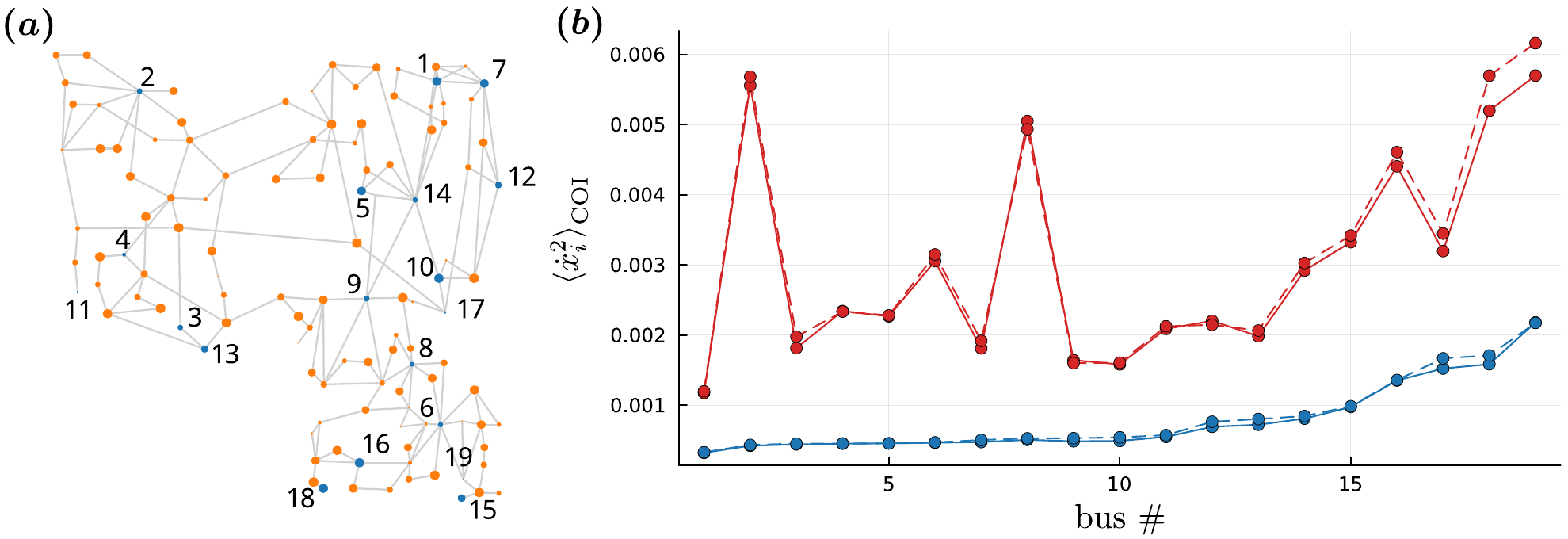}
    \caption{(a) Map of the IEEE 118 Testcase~\cite{118bus} system: slow and fast buses are displayed in blue and orange respectively, and the dots' sizes represent the standard deviation of the input noise $\sigma_i$. 
    (b) Variances $\big\langle \dot x_i^2\big\rangle_{\rm COI}$ of the frequency deviation from the center of inertia, obtained through numerical simulations (solid) for correlated noises $\bm \xi$ (red) and naive noises $\bm \eta_\mathcal{S}$ (blue). 
    They are compared with analytical expressions (dashed) of: the full expression (red) in~(\ref{eqfr}) and only the contribution of slow nodes (blue). 
    Ornstein-Uhlenbeck noises were used with typical correlation times $\tau_i=0.1$  and standard deviations $\sigma_i \in [0, 0.01]$.}
    \label{fig:var}
\end{figure*}

In~(\ref{eqfr}), we set the standard deviation of the ambient noise in the slow and fast components to $\bm \sigma_{\mathcal{S}}$ and $\bm \sigma_{\mathcal{F}}$ respectively. We also set distinct homogeneous correlation times for the noise in each component as $\tau_i=\tau_\mathcal{S}$ for $i\in \mathcal{S}$ and $\tau_i=\tau_\mathcal{F}$ for $i\in \mathcal{F}$. While the contribution to the variance from the additive noise in the slow component is essentially given by the position of the buses on the slowest eigenmodes, the effect of the noise coming from the fast component involves combinations of eigenmodes. The precise combination depends on the effective reduced dynamics through $\Gamma_{\alpha\beta}$.

\begin{figure*}[h!]
    \centering
    \includegraphics[width=0.99\linewidth]{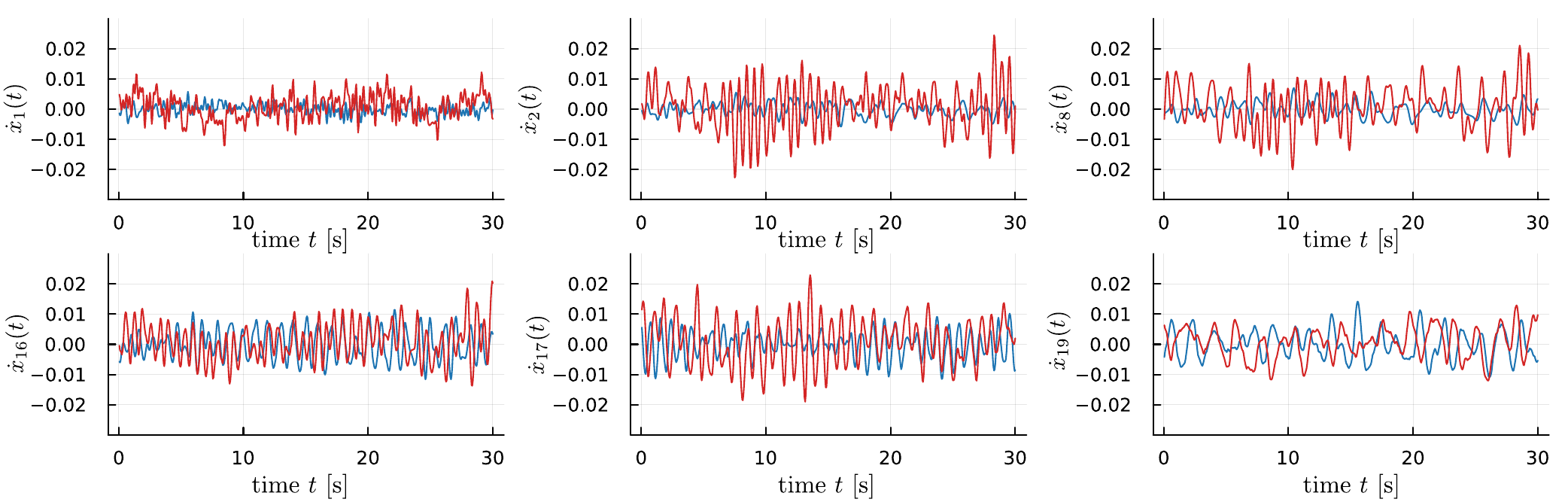}
    \caption{Temporal evolution of frequencies $\dot x_i(t)$ for a selection of 6 generator buses. 
    Noises are simulated as the naive white noise $\bm{\eta}_{\cal S}$ (blue) and the corrected correlated noise $\bm{\xi}$. 
    These time series were generated with the same $\sigma_i$ as in Fig.~\ref{fig:var}. 
    The approach presented in this work indeed produces frequency patterns quite different from the naive approach.}
    \label{fig:xdot}
\end{figure*}

\subsection{Specific Grid Topologies}\label{ssec:topol}

To get an intuition of the underlying mathematics that lead (or not) to a discrepancy between the dynamics in the reduced and full systems, we consider two idealized cases. Of course, such idealized systems are not really representative of actual implementations, but they are analytically tractable and provide insight into the Kron reduction's behavior.

Let us consider the case of a star graph where a load node is surrounded by $N_{\cal S}$ generators, as shown in Fig.~\ref{fignet} (a). 
In this case, one readily verifies that, provided phase differences are small,
\begin{align}
 \!\!\!   \JSS &\approx -\Id_{N_{\cal S}}\, , & \JFF &\approx -N_{\cal S}\, , & \JSF &= \JFS^\top \approx \bm{1}_{N_{\cal S}}\, ,
\end{align}
which yields $\Gamma_{\alpha\beta} \approx 0$ independently of $\alpha$ and $\beta$\,, where we assumed homogeneous standard deviation for the noise in the fast component. In such a case, the error induced by the reduction is rather small.

On the contrary, when a single generator node is connected to $N_{\cal F}$ loads, as shown in Fig.~\ref{fignet} (b), one gets
\begin{align}
 \!\!\!\!\!      \JSS &\approx -N_{\cal F}\, , & \JFF &\approx -\Id_{N_{\cal F}}\, , & \JFS &= \JSF^\top \approx \bm{1}_{N_{\cal F}}\, .
\end{align}
Here, the coefficient $\Gamma_{\alpha\beta} = \sigma^2 N_{\cal F}$ (if $\alpha=\beta=1$, otherwise $\Gamma_{\alpha\beta} \approx 0$) scales linearly with the number of loads. The error in such a reduced system is typically large.

While the two cases depicted above are quite extreme and somewhat unrealistic, they provide interesting insights into situations where reducing a system may lead to more or less accurate estimates of the frequency variance. 
To put it shortly, a reduced node connected to many non-reduced ones will not lead to a large error, while the opposite is expected to induce a large error in the estimate. Power grids structures are typically closer to the latter situation also including connections between the reduced buses.  One therefore expects a non-negligible impact of the reduced buses on the synchronous generator dynamics. We numerically confirm this conjecture in the next section.

\section{Numerical Illustration}\label{sec:sim}
Here we compare the frequency deviations obtained with the full noise term $\bm \xi$ in (\ref{eq5}) and its naive approximation using uncorrelated noise, i.e., taking into account only ${\bm \eta}_{\mathcal{S}}$\,. 
To illustrate our findings, we use the standard IEEE 118 Testcase~\cite{118bus}. 
An OPF was used to obtain the initial steady-state configuration, i.e., phase differences and generators' outputs. 
The original test case was augmented with dynamical parameters deduced from generators' and loads' characteristics, similarly to~\cite{pagnier2019inertia}. 
For the sake of comparing numerical results with the theory, dynamical parameters were uniformized, leading to $\overline{d} = 0.0005~\text{[s]} \ll \underline{d}= 0.05~\text{[s]}$, and $\overline{m}= 0.002~\text{[s$^2$]} \ll \underline{m} = 0.2~\text{[s$^2$]}$.

Figure~\ref{fig:var} compares the stochastic differential equation proposed in this work with the current standard in the literature, which neglects the contributions from the noise on the reduced buses. 
Slow buses were ordered from the smallest variance according to the naive approach. 
Fig.~\ref{fig:var} (a) shows the system’s map and the amount of noise present at each bus. 
Notice that the system response to the external inputs is scattered thought it and that the variance is not directly related to the local noise and thus cannot be deduced solely from it. 
Fig.~\ref{fig:var} (b) shows the theoretical and empirical variance with (red) and without (blue) the correction term of \eqref{eqfr}. 
Comparing the red and blue curve, we see that both measures generally follow the same trend, which depends on the structure of the system. 
However, for buses that are in the bulk of slow buses, the naive approach can yield results that are not only quantitatively different but even qualitatively incorrect, leading to a wrong assessment of the system’s vulnerabilities. 
As a striking example, bus \#2 jumps from the penultimate position to the forefront.

Following the discussion of Sec.~\ref{ssec:topol} about specific grid topologies, we observe that the two buses with the largest discrepancies, i.e., \#2 and \#8, are in a starlike load configuration. 
Bus \#2 is isolated from the other slow buses and is thus affected by the noise of a large number of fast buses. 
On the other hand, bus \#14 is also in a starlike configuration, but in its case, most of its neighbors are other slow buses and are therefore less impacted.
We illustrate the striking discrepancies between naive and corrected trajectories in Fig.~\ref{fig:xdot}. 

\begin{figure}[h!]
    \centering
    \includegraphics[width=0.99\columnwidth]{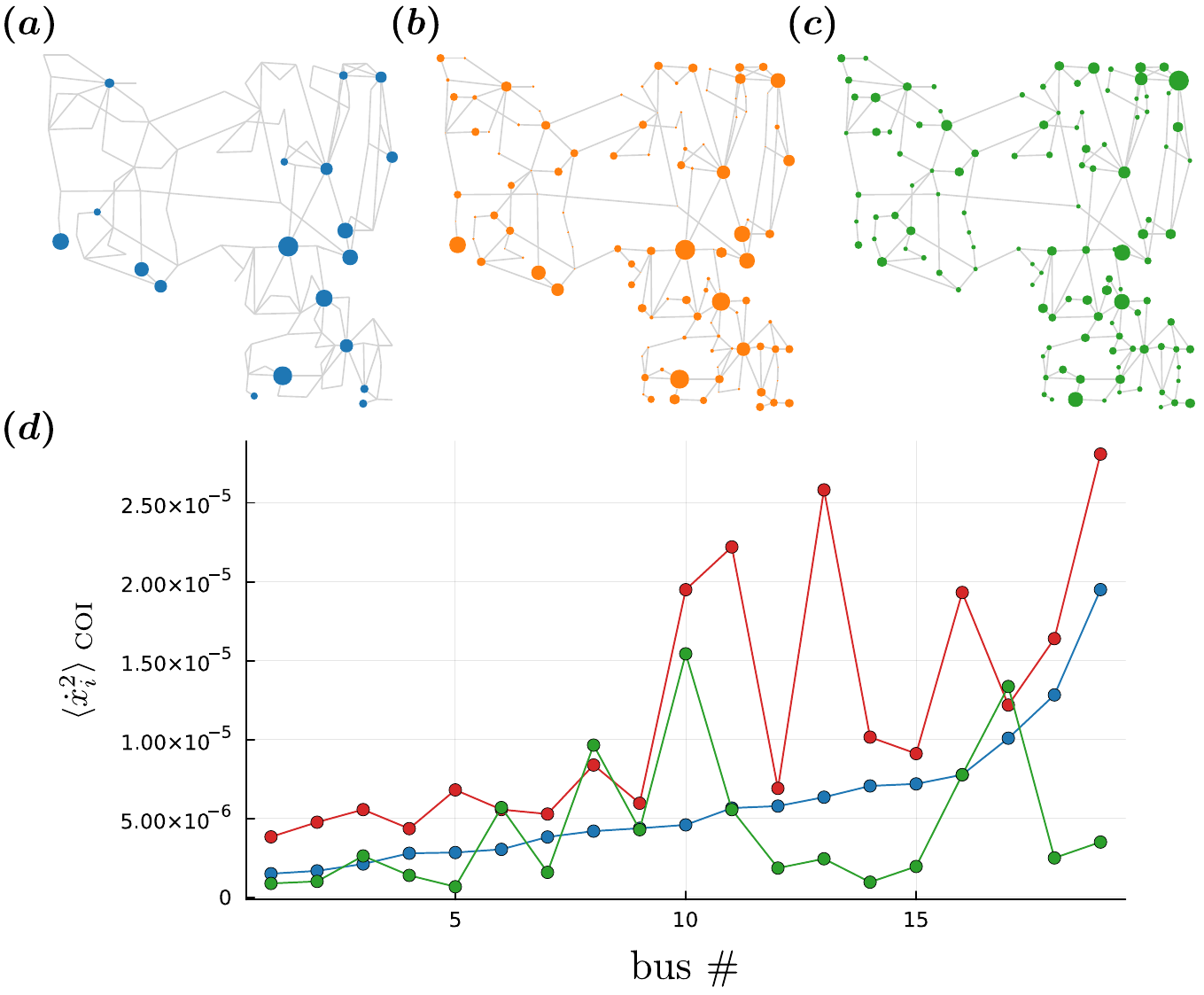}
    \caption{(a) Inertia coefficients $m_i$. (b) Damping coefficients $d_i$. (c) Noise standard deviations $\sigma_i$. (d) Variances $\langle \dot x_i^2\rangle_{\rm COI}$ computed with the reduced model consisting of inertial buses with the correlated noise $\bm \xi$ (red) and the naive approximation $\bm \eta_{\mathcal{S}}$ (blue), compared with those obtained with the full structure-preserving model (green).}
    \label{fig:hetero}
\end{figure}

Up to this point, we have assumed that inertia and damping were homogeneous. 
This assumption allowed us to describe the system dynamics in terms of the eigenvalues and eigenvectors of the Laplacian matrix and to derive a closed-form solution thereof, as well as for a resilience measure. 
If this assumption is not met, these eigenmodes are not independent and no direct interpretation can be made. 
However, the Mori-Zwanzig reduction is still applicable and \eqref{eq5} still holds. 
In Fig.~\ref{fig:hetero}, we compare both reduced models, with and without correlated noise, to the structure-preserving model~\cite{Ber81}. 
None of the models is quite able to correctly assess the variances, which is not surprising given that the system lies beyond the assumptions made here. 
Nonetheless, the model developed in this work seems to qualitatively better capture the system vulnerabilities and is therefore probably more suitable for use.

\section{Conclusion}\label{sec:conclusion}
We have shown that the Kron reduction, commonly used in power grid modeling, can lead to misleading estimates of the voltage dynamics. 
Typically, when the system is subject to noise, the variance of the voltage frequency can be significantly underestimated if the impact of reduced buses is neglected. 
Such underestimation can be of major importance when assessing the robustness of the grid against disturbances and could lead to safety issues. 
Accordingly, we recommend taking the reduced buses into account in simulations.

To address these limitations of the Kron reduction, we propose a method to incorporate the impact of the fast buses on the slow buses. 
Specifically, we use the Mori-Zwanzig formalism to derive the equation of motion for the slow buses as a function of a timescale parameter $\epsilon$\,. 
This approach elucidates how both the dynamics and the inputs at the slow buses affect the evolution of the slow ones. 
Expanding at the leading order in $\epsilon$\,, we recover for the reduced dynamics the usual Kron reduced Jacobian matrix. 
However, even if the initial inputs are uncorrelated, the effective noise acting on the non-reduced buses is correlated. 
In general, the correlation of the noise cannot by neglected. We show for the IEEE 118 Testcase that neglecting the correlation and considering only the noise coming from the non-reduced buses leads to a misleading assessment of the vulnerabilities. 
Therefore, to accurately evaluate the grid response, one should not solely use uncorrelated noise at the non-reduced buses but also take into account the contribution of external inputs at the reduced nodes. 
Importantly, our framework gives a way to go beyond the leading order in $\epsilon$ and calculate higher order correction.

Future work should consider a more accurate reduction for the full structure preserving model. The latter could be achieved by considering corrections beyond the leading order in $\epsilon$\,. 

\appendix

The numerical simulations were performed with \textit{DifferentialEquations.jl}~\cite{rackauckas2017differentialequations}, a Julia package that gathers well-optimized solvers often benchmarked as some of the fastest available implementations. Specifically, we used ImplicitEM, an order 0.5 Ito drift-implicit method with the Trapezoid method on the drift term. It is particularly suitable for stiff stochastic differential equations, which power system dynamics with realistic parameters tend to behave like. To ensure good agreement between analytical and numerical results, an upper bound had to be set on the time steps. Relaxing this constraint renders the scheme significantly faster, but at the expense of accuracy.


\begin{thebibliography}{10}
\providecommand{\url}[1]{#1}
\csname url@rmstyle\endcsname
\providecommand{\newblock}{\relax}
\providecommand{\bibinfo}[2]{#2}
\providecommand\BIBentrySTDinterwordspacing{\spaceskip=0pt\relax}
\providecommand\BIBentryALTinterwordstretchfactor{4}
\providecommand\BIBentryALTinterwordspacing{\spaceskip=\fontdimen2\font plus
\BIBentryALTinterwordstretchfactor\fontdimen3\font minus
  \fontdimen4\font\relax}
\providecommand\BIBforeignlanguage[2]{{%
\expandafter\ifx\csname l@#1\endcsname\relax
\typeout{** WARNING: IEEEtran.bst: No hyphenation pattern has been}%
\typeout{** loaded for the language `#1'. Using the pattern for}%
\typeout{** the default language instead.}%
\else
\language=\csname l@#1\endcsname
\fi
#2}}

\bibitem{brummitt2013transdisciplinary}
C.~D. Brummitt, P.~D.~H. Hines, I.~Dobson, C.~Moore, and R.~M. D'Souza,
  ``Transdisciplinary electric power grid science,'' \emph{Proc. Nat. Acad.
  Sci. USA}, vol. 110, no.~30, pp. 12\,159--12\,159, 2013.

\bibitem{smith2022effect}
O.~Smith, O.~Cattell, E.~Farcot, R.~D. O'Dea, and K.~I. Hopcraft, ``The effect
  of renewable energy incorporation on power grid stability and resilience,''
  \emph{Sci. Adv.}, vol.~8, no.~9, p. eabj6734, 2022.

\bibitem{bamieh2012coherence}
B.~Bamieh, M.~R. Jovanovic, P.~Mitra, and S.~Patterson, ``Coherence in
  {{Large-Scale Networks}}: {{Dimension-Dependent Limitations}} of {{Local
  Feedback}},'' \emph{IEEE Trans. Autom. Control}, vol.~57, no.~9, pp.
  2235--2249, 2012.

\bibitem{siami2014systemic}
M.~Siami and N.~Motee, ``Systemic measures for performance and robustness of
  large-scale interconnected dynamical networks,'' in \emph{53rd {{IEEE
  CDC}}}.\hskip 1em plus 0.5em minus 0.4em\relax Los Angeles, CA, USA: IEEE,
  2014, pp. 5119--5124.

\bibitem{tegling2015price}
E.~Tegling, B.~Bamieh, and D.~F. Gayme, ``The {{Price}} of {{Synchrony}}:
  {{Evaluating}} the {{Resistive Losses}} in {{Synchronizing Power
  Networks}},'' \emph{IEEE Trans. Control Netw. Syst.}, vol.~2, no.~3, pp.
  254--266, 2015.

\bibitem{grunberg2018performance}
T.~W. Grunberg and D.~F. Gayme, ``Performance {{Measures}} for {{Linear
  Oscillator Networks Over Arbitrary Graphs}},'' \emph{IEEE Trans. Control
  Netw. Syst.}, vol.~5, no.~1, pp. 456--468, 2018.

\bibitem{kron1939tensor}
G.~Kron, \emph{Tensor analysis of networks}.\hskip 1em plus 0.5em minus
  0.4em\relax J. Wiley \& Sons New York, 1939.

\bibitem{dorfler2013kron}
F.~D{\"o}rfler and F.~Bullo, ``Kron {{Reduction}} of {{Graphs With
  Applications}} to {{Electrical Networks}},'' \emph{IEEE Trans. Circuits Syst.
  I}, vol.~60, no.~1, pp. 150--163, 2013.

\bibitem{Hor94}
R.~A. Horn and C.~R. Johnson, \emph{Matrix Analysis}.\hskip 1em plus 0.5em
  minus 0.4em\relax New York: Cambridge University Press, 1994.

\bibitem{mori1965transport}
H.~Mori, ``Transport, {Collective} {Motion}, and {Brownian} {Motion},''
  \emph{Progr. Theor. Phys.}, vol.~33, no.~3, pp. 423--455, 1965.

\bibitem{zwanzig1973nonlinear}
R.~Zwanzig, ``Nonlinear generalized {Langevin} equations,'' \emph{J. Stat.
  Phys.}, vol.~9, no.~3, pp. 215--220, 1973.

\bibitem{Tyl24}
M.~Tyloo, ``Resilience of the slow component in timescale-separated
  synchronized oscillators,'' \emph{Frontiers Netw. Physiol.}, vol.~4, 2024.

\bibitem{Ber81}
A.~R. Bergen and D.~J. Hill, ``A structure preserving model for power system
  stability analysis,'' \emph{IEEE Trans. Power App. Syst.}, no.~1, pp. 25--35,
  1981.

\bibitem{Kun22}
P.~Kundur and O.~P. Malik, \emph{Power system stability and control}, second
  edition~ed.\hskip 1em plus 0.5em minus 0.4em\relax New York: McGraw Hill
  Education, 2022.

\bibitem{ronellenfitsch2018optimal}
H.~Ronellenfitsch, J.~Dunkel, and M.~Wilczek, ``Optimal noise-canceling
  networks,'' \emph{Phys. Rev. Lett.}, vol. 121, no.~20, p. 208301, 2018.

\bibitem{hindes2023stability}
J.~Hindes, I.~B. Schwartz, and M.~Tyloo, ``Stability of kuramoto networks
  subject to large and small fluctuations from heterogeneous and spatially
  correlated noise,'' \emph{Chaos}, vol.~33, no.~11, 2023.

\bibitem{118bus}
R.~Christie, ``{IEEE 118 Tescase},''
  \url{http://labs.ece.uw.edu/pstca/pf118/pg_tca118bus.htm}, 1993.

\bibitem{pagnier2019inertia}
L.~Pagnier and P.~Jacquod, ``Inertia location and slow network modes determine
  disturbance propagation in large-scale power grids,'' \emph{PloS one},
  vol.~14, no.~3, p. e0213550, 2019.

\bibitem{rackauckas2017differentialequations}
C.~Rackauckas and Q.~Nie, ``Differential{E}quations.jl--a performant and
  feature-rich ecosystem for solving differential equations in {J}ulia,''
  \emph{J. Open Research Software}, vol.~5, no.~1, 2017.

\end{thebibliography}
\end{document}